\RequirePackage{fix-cm}
\documentclass[smallextended]{svjour3}       
\smartqed  
\usepackage{graphicx}
\usepackage{subfigure}
\usepackage{amsfonts}

\newcommand{\Ii}{\mathbf{1}}
\spnewtheorem{algorithm}{Algorithm}[section]{\bf}{\rm}
\spnewtheorem{condition}{Condition}[section]{\bf}{\rm}
\begin{document}
\title{
WEIGHTED MAXIMA AND SUMS OF NON-STATIONARY RANDOM LENGTH SEQUENCES IN HEAVY-TAILED MODELS
} 

\author{Natalia M. Markovich
}
\institute{N. M. Markovich \at
              V.A.Trapeznikov Institute of Control Sciences,
                Russian Academy of Sciences,
                Profsoyuznaya 65, 
                117997 Moscow, Russia
                \\
              Tel.: +007-495-3348820\\
              Fax: +007-495-3342016\\
              \email{nat.markovich@gmail.com}           
           }

\date{Received: date / Accepted: date}

\maketitle

\begin{abstract}
The sums and maxima of weighted non-stationary random length sequences of regularly varying random variables   may have the same tail and extremal indices, Markovich and Rodionov (2020).
 The main constraints  are that there exists a unique series in a scheme of series with the minimum tail index, the tail of the term number is lighter than the tail of the terms  and the weights are positive constants. These assumptions are changed here: a bounded
random number of series is allowed to have the minimum tail index, the tail of the term number may be heavier than the tail of the terms  and the weights may be real-valued. Then we derive the tail and extremal indices of the weighted non-stationary random length sequences under the new assumptions.
\end{abstract}
\section{Introduction}\label{Sec1}
\par
 Random length sequences and distribution tails of their sums and maxima attract the interest of many researchers due to numerous applications including  queues, branching processes and random networks, see Asmussen and Foss (2018), Jessen and  Mikosch (2006), Goldaeva and Lebedev (2018), Lebedev (2015a, 2015b), Markovich and Rodionov (2020), Olvera-Cravioto (2012), Robert and Segers (2008), Tillier and Wintenberger (2018).
\par
Considering weighted sums and maxima we aim to extend the results obtained in  Markovich and Rodionov (2020). As in the latter paper, we deal with a doubly-indexed array $\{Y_{n,i}: n,i\ge 1\}$ of nonnegative random variables (r.v.s) in which the "row index" $n$ corresponds to the time, and the "column index" $i$  enumerates the series. On the same probability space, the existence of a sequence of non-negative integer-valued r.v.s $\{N_n: n\ge 1\}$ is assumed.
Let $\{Y_{n,i}: n\ge 1\}$ be a strict-sense stationary sequence  with the extremal index $\theta_i$ having a regularly varying tail 
\begin{eqnarray}\label{11a} P\{Y_{n,i}>x\}&=&\ell_i(x)x^{-k_i}\end{eqnarray}
with tail index $k_i>0$ and  a slowly varying function $\ell_i(x)$ for any $i\ge 1$. There are no assumptions on the dependence structure in $i$.
\begin{definition}\label{Def1} A strictly stationary sequence  $\{Y_n\}_{n\ge 1}$ with distribution function $F(x)$ and $M_n=\bigvee_{j=1}^{n}Y_j=\max_j Y_j $ is said to have the extremal index $\theta\in[0,1]$ if
for each $0<\tau <\infty$ there exists a sequence of real numbers $u_n=u_n(\tau)$ such that
\begin{eqnarray}\label{1}&&\lim_{n\to\infty}n(1-F(u_n))=\tau,\end{eqnarray}
\begin{equation}\label{2}\lim_{n\to\infty}P\{M_n\le u_n\}=e^{-\tau\theta}\end{equation}
hold (Leadbetter et al. (1983), p.63).
\end{definition}
I.i.d. r.v.s $\{Y_n\}$ give $\theta=1$. The converse may be incorrect.
An extremal index that is close to zero implies a kind of a strong local dependence.
\par
In Markovich and Rodionov (2020), the weighted sums and maxima
\begin{eqnarray}\label{3}
&&Y_{n}^{*}(z, N_n) =\max(z_1Y_{n,1},...,z_{N_n}Y_{n,N_n}),
\\
&&
~~Y_{n}(z, N_n)=z_1Y_{n,1}+...+z_{N_n}Y_{n,N_n}\nonumber\end{eqnarray}
for positive constants $z_1, z_2,...$ and regularly varying $\{Y_{n,i}: n\ge 1\}$ were considered. A similar result was obtained in Goldaeva (2013) for the same weights  and for random sequences of a fixed length $l\ge 1$ and when $\{Y_{n,i}: n\ge 1\}$ have a power-type tail, i.e. $P(Y_{n,i}>x) \sim c^{(i)} x^{-k_i}$ as $x\to\infty$,\footnote{The symbol $\sim$ means asymptotically equal to or $f(x)\sim g(x)$ $\Leftrightarrow$ $f(x)/g(x)\rightarrow 1$ as $x\rightarrow a$, $x\in M$ where the functions  $f(x)$ and $g(x)$ are defined on some set $M$ and $a$ is a limit point of $M$.} where $c^{(i)}$ is a 
positive constant for each fixed $i$.
\\
Let us recall Theorem 4 derived in Markovich and Rodionov (2020) and$\{Y_{n,i}: n\ge 1\}$ related to 
sums and maxima (\ref{3}) in the following Theorem \ref{T1}. It is assumed that the "column" sequences $\{Y_{n,i}: i\ge 1\}$
have stationary distribution tails (\ref{11a})
 in $n$ with positive tail indices $\{k_1, k_2,...\}$ and
extremal indices $\{\theta_1, \theta_{2},...\}$ for each fixed $i$. Here 
$\{\ell_i(x)\}$ are restricted by the 
condition: for all
$A>1$, $\delta>0$
there exists $x_0(A, \delta)$ such that for all $i\geq 1$
\begin{eqnarray} \ell_i(x)\leq A x^\delta,\ \  x>x_0(A, \delta) \label{uniform}
\end{eqnarray}
holds. 
$N_n$ has a regularly varying distribution
 with the tail index $\alpha>0,$ that is
 \begin{eqnarray}\label{15a} &&P(N_n>x) = x^{-\alpha} \tilde{\ell}_n(x).\end{eqnarray}
 There is a minimum tail index $k_1$ and
$k:= \lim_{n\to\infty} \inf_{2\leq i\leq l_n} k_i,$ 
 \begin{eqnarray}\label{27}&& l_n=[n^\chi],
  \qquad
   \end{eqnarray}
  and  $\chi$ satisfies
    \begin{equation}
0<\chi<\chi_0, \qquad\chi_0 = \frac{k-k_1}{k_1(k+1)}.
\label{chi}
\end{equation}
An arbitrary dependence structure  between  $\{Y_{n,i}\}$ and  $\{N_n\}$  is allowed. The tail of  $N_n$ does not dominate the tail of the most heavy-tailed term $Y_{n,1}$. Let  $\ell(x)$ be such that $\ell_1(x) = \ell^{-k_1}(x)$ and $\ell^{\sharp}(x)$ be the de Brujin conjugate of $\ell(x)$ (see Bingham et al. (1987)). Let
\begin{eqnarray}\label{23}u_n&=&yn^{1/k_1}\ell_1^{\sharp}(n), ~~y>0,\end{eqnarray} where we denote $\ell_1^{\sharp}(x) = \ell^{\sharp}(x^{1/k_1})$ and the positive weights $\{z_i\}$ are
bounded. 
\begin{theorem}\label{T2} (Markovich and Rodionov 2020)
Let $k_1<k$, (\ref{uniform}), (\ref{27})
and (\ref{chi}) hold. Then the sequences $Y_n^*(z,l_n)$ and $Y_n(z,l_n)$ have the same tail index
$k_1$ and the same extremal index $\theta_1$.
\end{theorem}
Theorem \ref{T1} follows by Theorem \ref{T2}.
 \begin{theorem}\label{T1}  (Markovich and Rodionov 2020)
Let the sets of slowly varying functions $\{\tilde{\ell}_n(x)\}_{n\geq 1}$ in (\ref{15a}) and $\{\ell_i(x)\}_{i\geq 1}$ in (\ref{11a}) satisfy the condition (\ref{uniform}).
Suppose that $k_1<k$ and
\begin{eqnarray}\label{4a}P\{N_n>l_n\}&=&o\left(P\{Y_{n,1}>u_n\}\right), ~~n\to\infty\end{eqnarray}
hold, where the sequence $l_n$ satisfies (\ref{27}) and (\ref{chi}).
Then the sequences $Y_{n}^*(z,N_n)$ and $Y_{n}(z,N_n)$
have the same tail index
$k_1$ and the same extremal index $\theta_1$.
\end{theorem}
\begin{remark} By the proof of the latter theorems it follows that if the "column" series with a minimum tail index $k_1$ is unique, then the tail distributions of the sequences $Y_{n}^*(z,N_n)$ and $Y_{n}(z,N_n)$ are asymptotically equivalent to its distribution tail. This is enough for $Y_{n}^*(z,N_n)$ and $Y_{n}(z,N_n)$ to have an extremal index $\theta_1$, despite  both sequences may be non-stationary distributed due to an arbitrary dependence between the "column" sequences.
\end{remark}
Our objectives are 
to revise Theorem \ref{T1} for the case when a random number of "column" series may have the minimum tail index both for positive and real-valued constant weights in (\ref{3}). Moreover, we observe how the results may change under  the opposite assumption 
\begin{eqnarray}\label{4b}P\{Y_{n,1}>u_n\}&=&o\left(P\{N_n>l_n\}\right), ~~n\to\infty\end{eqnarray}
instead of (\ref{4a}).
\\
We use the following notations
\begin{eqnarray*}M_n^{(i)}&=& \max\{Y_{1,i}, Y_{2,i},...,Y_{n,i}\}, ~i\in\{1,..,l_n\},
\\
M_n(z,l_n)&=&\max\{Y_1(z,l_n),Y_2(z,l_n),...,Y_n(z,l_n)\},
\\
M_n^*(z,l_n)&=&\max\{Y_1^*(z,l_n),...,Y_n^*(z,l_n)\}=\max\{z_1M_n^{(1)},...,z_{l_n}M_n^{(l_n)}\},~~n\ge 1.
\end{eqnarray*}
The paper is organized as follows. The revision of Theorems \ref{T2} and \ref{T1} by Theorems \ref{T3} and \ref{T4} for fixed and  random numbers of the most heavy-tailed "column" series and positive weights in (\ref{3}) is given in Section \ref{Sec2.1}. The same revision  for real-valued   weights is given in Section \ref{Sec2.2}. The revision of Theorem \ref{T1} for a heavy-tailed number $N_n$ of light-tailed terms is stated in Section \ref{Sec2.3}.  
 The proofs are given in Section \ref{Sec3}.

\section{Main Results}\label{Sec2}
\subsection{Revision of Theorem \ref{T1} for positive constant weights}\label{Sec2.1}

We  revise Theorem \ref{T1} by Theorem \ref{T4}  allowing  a random bounded number $d\ge 1$ of series to have a minimum tail index.
To this end, we extend Theorem \ref{T2} by Theorem \ref{T3}. 
Theorem \ref{T2} covers the case $d=1$. We assume in Theorem \ref{T3} 
that $d>1$ is fixed 
and $k_i=k_1$, $i\in\{1,..., d\}$,  $1\le d\le l_n-1$,
$k_1<k=k_{d+1}$, where
\begin{eqnarray}\label{5}k &:=& \lim_{n\to\infty} \inf_{d+1\leq i\leq l_n} k_i\end{eqnarray}
holds.
We introduce the following 
conditions:
\begin{enumerate}
\item[(A1)] The stationary sequences $\{Y_{n,i}\}_{n\ge 1}$, $i\in\{1,...,d\}$ are mutually independent, and  independent of the sequences
$\{Y_{n,i}\}_{n\ge 1}$, $i\in\{d+1,...,l_n\}$. 
\\
\item[(A2)] Assume $\{Y_{n,i}\}_{n\ge 1}$, $i\in\{1,...,d\}$ satisfy the following conditions as $x\to\infty$
\begin{eqnarray}\label{6b}\frac{P\{Y_{n,i}>x\}}{x^{-k_1}\ell_1(x)}&\rightarrow& c_i, ~~i\in\{1,...,d\},
\end{eqnarray}
for some non-negative numbers $c_i$,
\begin{eqnarray}\label{6c}\frac{P\{Y_{n,i}>x, Y_{n,j}>x\}}{x^{-k_1}\ell_1(x)}&\rightarrow& 0, ~~ i\neq j, ~~i,j\in\{1,...,d\}.
\end{eqnarray}
By Lemma 2.1 in Davis and Resnik (1996) in conditions (A2) it holds
\begin{eqnarray}\label{6}\frac{P\{\sum_{i=1}^dY_{n,i}>x\}}{x^{-k_1}\ell_1(x)}&\rightarrow& \sum_{i=1}^dc_i, ~~x\to\infty
\end{eqnarray}
for any $n\ge 1$.
\\
\item[(A3)] Assume that for each $n\ge 1$ there exists $i\in\{1,...,d\}$ such that 
\begin{eqnarray}\label{15}\!\!\!\!\!\!\!\!&& 
P\{\max_{1\le j\le d, j\neq i}(z_{j}Y_{n,j})>x, z_{i}Y_{n,i}\le x\}
=o(P\{z_iY_{n,i}>x\}),~~x\to\infty
\end{eqnarray}
holds.
\\
\item[(A4)] Assume that there exists $i\in\{1,...,d\}$ such that it holds
\begin{eqnarray}\label{11b}\!\!\!\!\!&& 
P\{\max_{1\le j\le d, j\neq i}(z_{j}M_n^{(j)})>u_n, z_{i}M_n^{(i)}\le u_n\}
=o(1),~~n\to\infty.
\end{eqnarray}
\end{enumerate}
\begin{example}\label{Exam3} Let $z_1Y_{n,1}\ge z_2Y_{n,2}\ge ... \ge z_dY_{n,d}$ for a given $n\ge 1$ hold. If $z_1Y_{n,1}\le x$ holds, then $z_iY_{n,i}\le x$ for all $i\in\{2,3,...,d\}$ and $P\{\max_{2\le j\le d}(z_{j}Y_{n,j})>x, z_{1}Y_{n,1}\le x\}=0$ follows.  Since for any "row" sequence $\{z_iY_{n,i}: i\ge 1\}$ there is a maximum term, the condition (\ref{15})
follows.
\end{example}
\begin{example}\label{Exam2} The assumption (\ref{11b}) is particularly valid for the following $d$ "column" sequences. Let elements of each  "column" sequence be sums 
 of the corresponding elements of all previous "columns", i.e. $Y_{n,i}=\sum_{j=1}^{i-1}Y_{n,j}$
 and $z_1\le z_2\le ...\le z_d$. Each of the $d$ "column" sequences has the tail index $k_1$ that follows from the  statement of Theorem \ref{T3}. Since $M_n^{(1)}=M_n^{(2)}< M_n^{(3)}<...<M_n^{(d)}$ holds, then  $\sum_{j=1}^{d-1} P\{z_{j}M_n^{(j)}>u_n,z_{j+1}M_n^{(j+1)}\le u_n,...,z_{d}M_n^{(d)}\le u_n\}=0$ or equivalently (\ref{11b}) for $i=d$ follow. In the same way, (\ref{11b}) is valid for all $d$ "column" sequences such that $M_n^{(1)}\le M_n^{(2)}\le M_n^{(3)}\le...\le M_n^{(d)}$ holds.
\end{example}
Let $u_n$ in Theorem \ref{T3} be as in (\ref{23}). 
\begin{theorem}\label{T3} Let
(\ref{uniform}) hold for all $d+1\leq i\leq l_n$, (\ref{27}), (\ref{chi}) be satisfied.
If, in addition, (A1) or (A2) holds, then $Y_n^*(z,l_n)$ and $Y_n(z,l_n)$ have  tail index $k_1$. If, instead of (A1) or (A2), (A3) holds, then $Y_n^*(z,l_n)$ has the same tail index.
\begin{enumerate}
\item\label{Item1} If 
(A1) holds, then
 $Y_n^*(z,l_n)$  and $Y_n(z,l_n)$ have 
the same extremal index
\begin{eqnarray}\label{5aaab}\theta(z) &= & \sum_{j=1}^d\theta_jz_j^{k_1}/\sum_{j=1}^dz_j^{k_1}.
\end{eqnarray}
\item\label{Item2} 
If (A4) 
holds, then $Y_n^*(z,l_n)$ has 
the  extremal index $\theta_i$. 
If, additionally to (A4),
(A2) holds, then $Y_n(z,l_n)$ has the same extremal index 
 as $Y_n^*(z,l_n)$.
 \end{enumerate}
\end{theorem}
\begin{remark}
An arbitrary dependence among elements of the $d$
"columns" series with the minimum tail index in Item 2  of Theorem \ref{T3} (for instance, elements of odd "row" sequences coincide and elements of the even rows are i.i.d.)  leads to 
the non-stationary sequences of maxima  $Y_n^*(z,l_n)$ and sums $Y_n(z,l_n)$.
However,
the extremal index of  $Y_n^*(z,l_n)$ 
exists in case of (\ref{11b}) 
due to the asymptotic equivalence of the distributions of $M_n^*(z,d)$ 
and $M_n^{(i)}$.
\end{remark}
Now we reformulate Theorem \ref{T1}.
\begin{theorem}\label{T4}
Let the sets of slowly varying functions $\{\ell_i(x)\}_{d+1\le i\le l_n}$ in (\ref{11a}) and $\{\tilde{\ell}_n(x)\}_{n\geq 1}$ in (\ref{15a}) satisfy the condition (\ref{uniform}) and (\ref{27}), (\ref{chi}), (\ref{4a}) hold. Assume that $d$ and $\{Y_{n,i}\}$ are  independent.
\begin{enumerate}
\item [(i)]
Let $d$ be a bounded discrete r.v.  such that
$1<d<d_n = \min(C, l_n)$, 
$C>1$ holds. 
\begin{enumerate}
\item 
%
    If (A1) or (A2) for any $d\in\{2,3,...,\lfloor d_n-1\rfloor\}$ 
holds and
$N_n$ and $\{Y_{n,i}\}$ are  independent, then $Y_{n}(z,N_n)$ and $Y^*_{n}(z,N_n)$ have the  tail index $k_1$. If, instead of (A1) and (A2), (A3) holds, then $Y_n^*(z,N_n)$ has the same tail index.
\item If (A4) where in (\ref{11b}) $d$ is replaced by $\lfloor d_n-1\rfloor$
holds, then $Y_n^*(z,N_n)$ has the  extremal index $\theta_i$. If, in addition, (A1) (or (A2)) for any 
$d\in\{2,3,...,\lfloor d_n-1\rfloor\}$ holds, then $Y_n(z,N_n)$ has the same extremal index.
\end{enumerate}
\item [(ii)] Suppose that $d>1$ is a bounded discrete r.v. equal to a positive integer  a.s.. Then all statements of Item (i) are fulfilled.
\end{enumerate}
\end{theorem}
\subsection{Revision of Theorem
 \ref{T1} for real-valued constant weights}\label{Sec2.2}

Let us consider (\ref{3}) for real-valued constant weights. The sums and maxima may be partitioned into negative and positive parts
\begin{eqnarray}\label{3b}
&&Y_{n}^*(z, N_n^+, N_n^-) =\max(Y_{n}^{*+}(z^+, N_n^+),Y_{n}^{*-}(z^-, N_n^-))=Y_{n}^{*+}(z^+, N_n^+),\nonumber
\\
&&Y_{n}^{**}(z, N_n^+, N_n^-) =\min(Y_{n}^{*+}(z^+, N_n^+),Y_{n}^{*-}(z^-, N_n^-))=Y_{n}^{*-}(z^-, N_n^-),\nonumber
\\
&&
~~Y_{n}(z, N_n^+, N_n^-)=Y_{n}^+(z^+, N_n^+)+Y_{n}^-(z^-, N_n^-),\end{eqnarray}
where the terms marked by '$+$' and '$-$' correspond to sums and maxima with positive and negative 
weights and the corresponding random numbers of terms $N_n^+$ and $N_n^-$. $N_n^{\pm}$ satisfies (\ref{15a}). For brevity we denote in the section   $Y_{n}^*(z, N_n^+, N_n^-)$, $Y_{n}^{**}(z, N_n^+, N_n^-)$ and $Y_{n}(z, N_n^+, N_n^-)$ as $Y_{n}^*(z, N_n)$, $Y_{n}^{**}(z, N_n)$ and $Y_{n}(z, N_n)$, respectively.
\\
We avoid to use only positive weights, the case derived in Markovich and Rodionov (2020), or only negative 
weights  since then
 $P\{Y_{n}^*(z, N_n)>u_n\}=P\{Y_{n}^*(z^-, N_n^-)>u_n\}=0$ and $P\{Y_{n}(z, N_n)>u_n\}=P\{Y_{n}^-(z^-, N_n^-)>u_n\}=0$ hold for positive $u_n$.
 \\
 We assume that the number of column subsequences $l_{n^+}$ and $l_{n^-}$ 
 of $l_n$ such that $l_{n^+}+l_{n^-}=l_n$, 
 $n^{\pm}\in \mathbb{N}$, $n^{\pm}\to\infty$ for positive and negative weights, respectively, satisfy conditions similar to (\ref{27}) and (\ref{chi}), i.e.
\begin{eqnarray}\label{27a}&& l_{n^{\pm}}=[{n^{\pm}}^{\chi^{\pm}}], 
  \qquad
   \end{eqnarray}
    and  $\chi^{\pm}$ satisfies
    \begin{equation}
0<\chi^{\pm}<\chi_0^{\pm}, \qquad\chi_0^{\pm} = \frac{k^{\pm}-k_1^{\pm}}{k_1^{\pm}(k^{\pm}+1)},
\label{27b}
\end{equation}
where
\begin{eqnarray}\label{27c}
&& k^+:= \lim_{n\to\infty} \inf_{d^++1\leq i\leq l_{n^+}} k_i^+,\qquad k^-:= \lim_{n\to\infty} \inf_{d^-+1\leq i\leq l_{n^-}} k_i^-,
\end{eqnarray}
Without loss of generality, we assume that the first $d^+$ and $d^-$, $1\le d^{\pm}\le l_{n^{\pm}}-1$, 
series in the corresponding subsequences have the minimum tail indices
$k_1^{+}$ and $k_1^{-}$, $k_1^{\pm}<k^{\pm}$, respectively. 
\\
By (\ref{3b}) all statements of Theorems \ref{T1} and \ref{T4} regarding the maximum $Y_{n}^*(z, N_n)$ (or $-Y_{n}^{**}(z, N_n)$) 
are fulfilled with a replacement of $k_1$ by the minimum tail index $k_1^+$ (or $k_1^-$) of the positive (or negative) weighted column sequences subject to the required assumptions.
\subsubsection{A unique "column" series with a minimum tail index}
Let there be a unique column with a positive  weight which has a minimum tail index $k_1^{+}$ and the  extremal index $\theta_1^{+}$, and similarly, there be a unique column with a
negative weight which has a minimum tail index $k_1^{-}$ and the  extremal index $\theta_1^{-}$.
We select $u_n$ similar as in Theorem \ref{T1}, i.e.  $u_n^+ = y n^{1/k_1^+} \ell^{\sharp}_1(n)$ or $u_n^- = y n^{1/k_1^-} \ell^{\sharp}_1(n)$
for $y>0$. Let us denote the column series with the minimum tail indices $k_1^{+}$ and $k_1^{-}$ as $Y_{n,1}^{+}$ and $Y_{n,1}^{-}$, respectively.
\begin{corollary}\label{Cor3}
Let the sets of slowly varying functions $\{\ell_i(x)\}_{i\geq 1}$ in (\ref{11a}) and $\{\tilde{\ell}_n(x)\}_{n\geq 1}$ in (\ref{15a}) satisfy the condition (\ref{uniform}) and $k_1^{\pm}<k^{\pm}$ hold.
Suppose that 
\begin{eqnarray}\label{4ab}P\{N_n^{\pm}>l_{n^{\pm}}\}
&=&o\left(P\{Y_{n,1}^{\pm}>u_n^{\pm}\} 
\right), ~~n\to\infty,\end{eqnarray}
holds, where the subsequences $l_{n^{\pm}}$ 
satisfy (\ref{27a}) and (\ref{27b}). 
\begin{enumerate}
\item
If $k_1^+\le k_1^-$ holds, then the sequences $Y_{n}^*(z, N_n)$ and $Y_{n}(z, N_n)$ have the same tail index $k_1^+$  and the same extremal index $\theta_1^+$ for 
  $u_n^+$, and their extremal indices do not exist for $u_n^-$.
\item If $k_1^-\le k_1^+$ holds, then  the sequences $-Y_{n}^{**}(z, N_n)$ and $-Y_{n}(z, N_n)$ have the same tail index $k_1^-$  and the same extremal index $\theta_1^-$ for
  $u_n^-$, and their extremal indices do not exist for $u_n^+$.
\end{enumerate}
\end{corollary}

\subsubsection{A random number 
of "column" series with a minimum tail index}
Let 
a random number  of column series have a minimum tail index.
\begin{corollary}\label{Cor4} Let the conditions of Theorem \ref{T4} with (\ref{4ab}) instead of (\ref{4a}) be fulfilled,   $d^+$ and $d^-$ be bounded discrete r.v.s such that $d^{\pm}<d_n^{\pm}=\min(C^{\pm},l_n^{\pm})$, $C^{\pm}>1$  
hold and $d^+$ and $d^-$ be independent of $\{Y_{n,i}\}$.
\begin{enumerate}
\item
If $k_1^+\le k_1^-$ holds, then all statements of Theorem \ref{T4}  are fulfilled for 
sequences $Y_{n}^*(z, N_n)$ and $Y_{n}(z, N_n)$ for $u_n^{+}$, 
and the
extremal index of $Y_{n}^*(z, N_n)$ and $Y_{n}(z, N_n)$ does not exist for $u_n^{-}$.
\item If $k_1^-\le k_1^+$ holds, then all statements of the previous case are fulfilled
by substituting $Y_{n}(z, N_n)$ by  $-Y_{n}(z, N_n)$, $Y_{n}^*(z, N_n)$ by   $-Y_{n}^{**}(z, N_n)$, $k_1^+$ by $k_1^-$ and $u_n^{\pm}$ by $u_n^{\mp}$ symmetrically.
 \end{enumerate}
\end{corollary}

\subsection{Revision of Theorem \ref{T1} for a heavy-tailed number of light-tailed terms}\label{Sec2.3}
The assumption (\ref{4a}) is crucial for the proof of Theorem \ref{T1}. We will replace the latter assumption with an opposite one and derive new statements for a unique column sequence and a random number of column sequences with a minimum tail index 
 for positive 
 weights.
\subsubsection{A unique "column" series with a minimum tail index}
\begin{theorem}\label{T5} Let the conditions of Theorem \ref{T1} be fulfilled apart of (\ref{4a}) and the sequence $l_n$ satisfy (\ref{27}) and (\ref{chi}).
\\
(i) Assume (\ref{4b})
and
\begin{eqnarray}\label{4d}%
\alpha\chi_0&>& 1+\frac{\alpha}{k_1}\delta^*,
\end{eqnarray}
for  $\delta^*>0$ hold, where 
$\alpha>0$ be the tail index of $N_n$.
Then the 
sequences $Y_{n}^*(z,N_n)$ and $Y_{n}(z,N_n)$ 
have the same tail index $k_1$ and
if \begin{eqnarray}\label{4e}%
\alpha\chi_0&>& 1+\frac{1}{2}(1-\alpha\chi),
\end{eqnarray}
holds, then the same  extremal index $\theta_1$.
\\
(ii) Assume
\begin{eqnarray}\label{4c}P\{N_n>l_n\} &=& P\{z_1Y_{n,1}>u_n\}(1+o(1)) 
\end{eqnarray}
 holds, and $N_n$ is regularly varying (\ref{15a}). 
  Then $Y_{n}^*(z,N_n)$  and $Y_{n}(z,N_n)$ 
have the same tail index $k_1$
and  the extremal index $\theta_1$.
\end{theorem}
\subsubsection{A random number of  "column" series with a minimum tail index}
\begin{theorem}\label{T6} Let the conditions of Theorem \ref{T4} be fulfilled, but (\ref{4a}) is substituted by (\ref{4b}) and (\ref{4d}). Then 
all statements of Theorem \ref{T4} regarding tail and extremal indices are fulfilled.
\end{theorem} 
\section{Proofs}\label{Sec3}
\subsection{Proof of Theorem \ref{T3}} 
The proof extends 
Theorem 3 in Markovich and Rodionov (2020).  We just indicate the modifications.
 The numeration of Theorem 3 in Markovich and Rodionov (2020) is preserved throughout the proof.
For brevity, we denote in this section $Y_n(z)=Y_n(z,l_n)$, $Y_n^*(z)=Y_n^*(z,l_n)$, $M_n(z)=M_n(z,l_n)$ and $M_n^*(z)=M_n^*(z,l_n)$.
\\
The right-hand side of (12) in Markovich and Rodionov (2020) can be rewritten as
\begin{eqnarray}\label{14}\!\!\!\!P\{Y_n(z)>u_n\}&\le & P\{\sum_{i=1}^{d}z_iY_{n,i}>u_n(1-\varepsilon)\}+\sum_{i=d+1}^{l_n}P\{z_iY_{n,i}>u_n\varepsilon_i\},
\end{eqnarray}
where  $\sum_{i=1}^{l_n}\varepsilon_i=1$ holds and $\{\varepsilon_i\}$  is a sequence of positive elements. Let us denote $\varepsilon=\sum_{i=d+1}^{l_n}\varepsilon_i$.
One may take $\varepsilon_i$, $i\in\{d+1,...,l_n\}$, in such a way
to satisfy $\varepsilon_i\to 0$ and $\varepsilon\to 0$ as $n\to\infty$. Choosing $\{\varepsilon_i=1/l_n^{\eta+1}\}$, $\eta>0$ as in Lemma 1 in  Markovich and Rodionov (2020) one can derive (13) in Markovich and Rodionov (2020) substituting $2$ by $d+1$, namely, it holds
\begin{eqnarray}\label{5a}&&\sum_{i=d+1}^{l_n}P\{z_iY_{n,i}>u_n\varepsilon_i\}=o(1/n), \qquad n\to\infty.\end{eqnarray}
 To prove the latter, we need to assume (\ref{uniform}) for $d+1\leq i\leq l_n$. 
For the selected $u_n$ we have
\begin{eqnarray}P\{z_1 Y_{n,1} > u_n\}  = (z_1/y)^{k_1} n^{-1} (1+o(1)), ~~n\to\infty.\label{quantile}\end{eqnarray}
We obtain
\begin{eqnarray}\label{4}&&P\{\sum_{i=1}^{d}z_iY_{n,i}>u_n(1-\varepsilon)\}
\le \sum_{i=1}^{d}P\{z_iY_{n,i}>u_n(1-\varepsilon)\varepsilon_i^*\}\nonumber
\\
&=& \frac{n^{-1}}{(y(1-\varepsilon))^{k_1}} \sum_{i=1}^d \left(\frac{z_i}{\varepsilon_i^*}\right)^{k_1}(1+o(1))=\left(\frac{z^*}{y(1-\varepsilon)}\right)^{k_1}n^{-1}(1+o(1)), 
\end{eqnarray}
where $\sum_{i=1}^d \varepsilon_i^*=1$ and $(z^*)^{k_1}= \sum_{i=1}^d \left(z_i/\varepsilon_i^*\right)^{k_1}$
hold.
By (\ref{14})-(\ref{4}) it follows
\begin{eqnarray}\label{8a}
&&\left(z_1/y\right)^{k_1}n^{-1}(1+o(1))=P\{z_1Y_{n,1}>u_n\}\le P\{Y_n^*(z)>u_n\} \nonumber
\\
&\le &P\{Y_n(z)>u_n\}\le  
\left(z^*/y(1-\varepsilon)\right)^{k_1}n^{-1}(1+o(1)) + o(1/n).
\end{eqnarray}
If 
the condition (A2) holds and replacing $x$ in (\ref{6}) by $u_n$, we get for any $m\ge 1$
\begin{eqnarray*}P\{\sum_{i=1}^dz_iY_{m,i}>u_n\}&=&n^{-1}\sum_{i=1}^d\left(\frac{z_i}{y}\right)^{k_1}c_i(1+o(1))
\end{eqnarray*}
since by (\ref{6b}) and (\ref{quantile})
\begin{eqnarray*}\frac{P\{z_iY_{n,i}>u_n\}}{P\{z_1Y_{n,1}>u_n\}}&=&\left(\frac{z_i}{z_1}\right)^{k_1}c_i(1+o(1)),~~n\to\infty,~~i=1,2,...,d.
\end{eqnarray*}
By (\ref{14}) and (\ref{5a}) this implies
\begin{eqnarray}\label{8b}\lim_{n\to\infty}nP\{Y_n(z)>u_n\}&=&\lim_{n\to\infty}nP\{Y_n^*(z)>u_n\}=\sum_{i=1}^d\left(\frac{z_i}{y}\right)^{k_1}c_i.
\end{eqnarray}
The latter is valid for  $Y_m^*(z)$ by (\ref{6b}), (\ref{6c}) and due to
\begin{eqnarray*}&&P\{\max (z_1Y_{m,1},z_2Y_{m,2})>u_n\}
  \\
  &=&P\{z_1Y_{m,1}>u_n\}+P\{z_2Y_{m,2}>u_n\}-P\{z_1Y_{m,1}>u_n,z_2Y_{m,2})>u_n\}
  \end{eqnarray*}
   (the case for general $d$ follows by induction)
 and since it holds
\begin{eqnarray*}&& P\{Y_n^*(z,d)>u_n\}\le P\{Y_n^*(z)>u_n\}\le P\{Y_n(z,d)>u_n\}.
\end{eqnarray*}
\paragraph{\textbf{Tail index}}
 Let us find the tail index
 of  $Y_m(z)$ and $Y_m^*(z)$  for all $m\ge 1$.
 Notice that the relation (\ref{14}) remains true when replacing $u_n$ by $x>0$. Similarly to (\ref{4}) it holds
  \begin{eqnarray*}\label{18}P\{\sum_{i=1}^{d}z_iY_{m,i}>x(1-\varepsilon)\}
 &\le &\frac{1}{(x(1-\varepsilon))^{k_1}} \sum_{i=1}^d \left(\frac{z_i}{\varepsilon_i^*}\right)^{k_1}\ell_i(x)(1+o(1))\nonumber
 \\
 &=&\frac{\ell^*(x)}{x^{k_1}}(1+o(1)),~~x\to\infty,
\end{eqnarray*}
where $\ell^*(x)=\sum_{i=1}^d \left(z_i/(\varepsilon_i^*(1-\varepsilon))\right)^{k_1}\ell_i(x)$ is a slowly varying function.
   Formula (16) 
  in  Markovich and Rodionov (2020) is 
  fulfilled by replacing $2$ by $d+1$ in the sum, namely,
  \begin{eqnarray}\label{18b}\sum_{i=d+1}^{l_m}P\{z_i Y_{n,i}>x\varepsilon_i\}&\le & O\left(x^{-k_1(1+\delta)}\right),~~x\to\infty,
  \end{eqnarray}
  where $\delta\in(0,(k-k_1)/(k+1))$.
  Then it follows that
 \begin{eqnarray}\label{13a}&&(z_1/x)^{k_1}\ell_1(x)(1+o(1))=P\{Y_{m,1}(z)>x/z_1\}\le P\{Y^*_m(z)>x\}\nonumber
 \\
 \!\!\! &\le &P\{Y_m(z)>x\}\le  P\{\sum_{i=1}^{d}z_iY_{m,i}>x(1-\varepsilon)\}+ \sum_{i=d+1}^{l_m}P\{z_iY_{n,i}>x\varepsilon_i\}\nonumber
 \\
 &\le & \ell^*(x)x^{-k_1}(1+o(1)).
 \end{eqnarray}
 By (\ref{13a}) it does not follow that $Y_m(z)$ and $Y_m^*(z)$ are  regularly varying.
 \paragraph{Condition (A1)} The tail index of $Y_n(z)$ and $Y_n^*(z)$ is equal to $k_1$ if the condition (A1) holds. Really,
it follows
\begin{eqnarray}\label{34a}P\{Y_n^*(z)>x\}&=& \sum_{i=1}^{d}\left(\frac{z_i}{x}\right)^{k_1}\ell_i(x)(1+o(1)), ~~ x\to\infty
\end{eqnarray}
by (\ref{8}) with replacement $u_n$ by $x>0$, and
\begin{eqnarray}\label{34b}\!\!\!\!\!P\{Y_n(z)>x\}&=& \sum_{i=1}^{d}P\{z_iY_i> x\}(1+o(1))=\sum_{i=1}^{d}\left(\frac{z_i}{x}\right)^{k_1}\ell_i(x)(1+o(1))
\end{eqnarray}
as $x\to\infty$ by Lemma 3.1 in Jessen and Mikosch (2006) due to independence of the first $d$ "column" sequences.
 \paragraph{Conditions (A2)} In conditions (A2) and due to (\ref{6}), (\ref{18b}), (\ref{13a}) we get
  \begin{eqnarray}\label{13b}&&P\{Y_m(z)>x\}= x^{-k_1}\ell_1(x)\sum_{i=1}^dz_i^{k_1}c_i(1+o(1))
 \end{eqnarray}
 since $P\{Y_m(z)>x\}\ge P\{\sum_{i=1}^{d}z_iY_{m,i}>x\}$.
  The same is valid for  $Y_m^*(z)$  by (\ref{6b}) and (\ref{6c}) in the same way as above.
 Hence, $Y_m(z)$ and $Y_m^*(z)$  have the same tail index $k_1$ for arbitrary (fixed) $m\ge 1$. 
\paragraph{Condition (A3)} We have for all $m\ge 2$
\begin{eqnarray*}\label{41}&& P\{Y_m^*(z)>x\}=P\{z_{1}Y_{m,1}>x\}+\sum_{i=2}^{l_m}P\{z_{i}Y_{m,i}>x,\max_{1\le j\le i-1}z_{j}Y_{m,j}\le x\}.\nonumber
  \end{eqnarray*}
  Let us denote
  \begin{eqnarray}\label{51}
  P\{z_{i}Y_{m,i}>x,\max_{1\le j\le i-1}z_{j}Y_{m,j}\le x\}&=& \xi_i.
  \end{eqnarray}
  We obtain
  \begin{eqnarray*}\label{32b}
  \sum_{i=2}^{l_m}\xi_i &=&\sum_{i=2}^{d}\xi_i+\sum_{i=d+1}^{l_m}\xi_i.
    \end{eqnarray*}
   By (\ref{18b})
   we have
   \begin{eqnarray*}&&
   \sum_{i=d+1}^{l_m}\xi_i\le \sum_{i=d+1}^{l_m}P\{z_{i}Y_{m,i}>x\}\le O\left(x^{-k_1(1+\delta)}\right),~x\to\infty.
   \end{eqnarray*}
  By (\ref{15}) it holds
  \begin{eqnarray*}
\sum_{i=2}^{d}\xi_i&=&\sum_{i=2}^{d}P\{z_iY_{m,i}>x,\max_{1\le j\le i-1}z_jY_{m,j}\le x\}
\\&=&P\{\max_{2\le j\le d}(z_{j}Y_{m,j})>x, z_{1}Y_{m,1}\le x\}=o(P\{z_{1}Y_{m,1}\le x\}),~~x\to\infty.
\end{eqnarray*}
Then it follows
\begin{eqnarray*} &&\!\!\!\!\!\!\!\!\!\!\!\!\!P\{Y^*_{m}(z))>x\}=P\{z_{1}Y_{m,1}>x\}(1+o(1))=(z_1/x)^{k_1}\ell_1(x)(1+o(1))
   \end{eqnarray*}
   and $Y^*_{m}(z)$ has tail index $k_1$.
\paragraph{\textbf{The extremal index for $d$ independent "column" sequences}}
Let us  assume  that the condition (A1) hold.
Due to independence we have
\begin{eqnarray*}
&&P\{M_n^*(z)\le u_n\} = \prod_{j=1}^{d}P\{z_j M_n^{(j)}\le u_n\}P\{z_{d+1}M_n^{(d+1)}\le u_n,...,z_{l_n}M_n^{(l_n)}\le u_n\}.
\end{eqnarray*}
Since the extremal index of the sequence $\{Y_{n,j}\}_{n \geq 1}$, $1\le j\le d$  
is assumed to be equal to $\theta_j,$ by (\ref{2}) and (\ref{quantile})
 we get
\begin{eqnarray*}\lim_{n\to\infty}\prod_{j=1}^{d}P\{z_jM_n^{(j)}\le u_n\}&=& \exp(-\sum_{j=1}^{d}\theta_j(z_j/y)^{k_1}). \end{eqnarray*}
By (19) and (21) in Markovich and Rodionov (2020) where $d+2$ 
replaces $2$ when indexing the  sums and by  (\ref{2}), 
 we obtain
\begin{eqnarray*}
\!\!\!\!\!&& P\{z_{d+1}M_n^{(d+1)}\le u_n,...,z_{l_n}M_n^{(l_n)}\le u_n\}= \!P\{z_{d+1} M_n^{(d+1)} \leq u_n\} (1+o(1))\to 1
\end{eqnarray*}
assuming
$k=k_{d+1}$,
since
\begin{eqnarray*}
&& nP\{z_{d+1}Y_{n,d+1}>u_n\}= n\cdot n^{-k_{d+1}/k_1}(z_{d+1}/y)^{k_{d+1}}\ell(n)(1+o(1))\to 0
\end{eqnarray*}
as $n\to\infty$ due to $k_{d+1}>k_1$ holds.
We obtain
\begin{eqnarray}\label{10}
&& \lim_{n\to\infty}nP\{Y_n^*(z)> u_n\}=\sum_{j=1}^d\left(\frac{z_{j}}{y}\right)^{k_1}
\end{eqnarray}
since it holds
\begin{eqnarray*}
&&P\{Y_n^*(z)> u_n\}=P\{\max(z_1Y_{n,1},...,z_{d}Y_{n,d},z_{d+1}Y_{n,d+1},...,z_{l_n}Y_{n,l_n})> u_n\}
\\
&=& 1-\prod_{i=1}^d(1-P\{z_iY_{n,i}> u_n\})(1-P\{\max(z_{d+1}Y_{n,d+1},...,z_{l_n}Y_{n,l_n})> u_n\}).
\end{eqnarray*}
By
\begin{eqnarray*}&&\prod_{i=1}^d(1-P\{z_iY_{n,i}> u_n\})=1-\sum_{i=1}^dP\{z_iY_{n,i}> u_n\}\\
&+&\sum_{i<j}P\{z_iY_{n,i}> u_n\}P\{z_jY_{n,j}> u_n\}-...-
(-1)^{d-1}\prod_{i=1}^dP\{z_iY_{n,i}> u_n\}
\\
&=& \left(1-\sum_{i=1}^dP\{z_iY_{n,i}> u_n\}\right)(1+o(1)),~~n\to\infty
\end{eqnarray*}
and by  (\ref{quantile}) it follows
\begin{eqnarray}\label{8}
&&P\{Y_n^*(z)> u_n\}\nonumber
\\
&=&\!\!\!\left(\sum_{i=1}^dP\{z_i Y_{n,i}> u_n\}+P\{\max(z_{d+1}Y_{n,d+1},...,z_{l_n}Y_{n,l_n})> u_n\}\right)\cdot (1+o(1))\nonumber
\\
&=& \sum_{i=1}^d\left(\frac{z_{i}}{y}\right)^{k_1}\cdot n^{-1}(1+o(1))+o(1/n) 
\end{eqnarray}
due to (9), (12), (13) and Lemma 1  in Markovich and Rodionov (2020) and $k=k_{d+1}$. 
By 
Lemma 3.1 in Jessen and Mikosch (2006)
it follows
\begin{eqnarray}\label{9}
&&\lim_{n\to\infty}nP\{Y_n(z)> u_n\}=\sum_{i=1}^d\left(\frac{z_{i}}{y}\right)^{k_1}.
\end{eqnarray}
By (\ref{10}) and
\begin{eqnarray*}\label{16}
&&P\{M_n^*(z)\le u_n\} =\exp(-\sum_{j=1}^{d}\theta_j(z_j/y)^{k_1}) (1+o(1))
\end{eqnarray*}
 we obtain that the extremal index of $Y_n^*(z)$ 
 is equal to (\ref{5aaab}).
In the same way as in the proof of Theorem 3 in Markovich and Rodionov (2020) one can show that $Y_n(z)$ has 
the same extremal index. 
 Really, we obtain
 \begin{eqnarray}\label{17}&&0\le P\{M_n^*(z)\le u_n\}-P\{M_n(z)\le u_n\}\nonumber
 \\
 &=&\sum_{j=1}^{n-1}P\{M_n^*(z)\le u_n, Y_j(z)>u_n, Y_{j+1}(z)\le u_n,...,Y_n(z)\le u_n\}\nonumber
 \\
 &+&P\{M_n^*(z)\le u_n, Y_n(z)> u_n\}\nonumber
 \\
 &\le & \sum_{j=1}^{n}P\{Y_j^*(z)\le u_n, Y_j(z)>u_n\}.
 \end{eqnarray}
 It follows
 \begin{eqnarray}\label{17a}\!\!\!\!\!\!\!\!\!\!P\{Y_j^*(z)\le u_n, Y_j(z)>u_n\}&=&P\{Y_j(z)>u_n\}-P\{Y_j^*(z)> u_n\}=o(n^{-1})
 \end{eqnarray}
 by
 (\ref{8})  and (\ref{9})
 for $j\le n$.
 \paragraph{\textbf{The extremal index for $d$ 
 dependent "column" sequences}}
We find 
the extremal index 
of the sequences $Y_n^*(z)$ and $Y_n(z)$ 
if there is the dependence of $\{Y_{n,i}\}_{n\ge 1}$ in $i$.  
\\
Let us rewrite (19)  in Markovich and Rodionov (2020) as
\begin{eqnarray}\label{12}
&&P\{M_n^*(z,d) 
\le u_n\}  -  \sum_{i=d+1}^{l_n}P\{z_i M_n^{(i)}> u_n\}\nonumber
\\
&\leq & P\{z_1M_n^{(1)}\le u_n,...,z_{l_n}M_n^{(l_n)}\le u_n\} = P\{M_n^*(z)\le u_n\} \nonumber
\\
&\leq & P\{M_n^*(z,d) 
\le u_n\}\leq P\{z_dM_n^{(d)}\le u_n\}.
\end{eqnarray}
By (21) in Markovich and Rodionov (2020) and assuming (\ref{uniform}), we have
\begin{eqnarray}\label{12a}\sum_{i=d+1}^{l_n}P\{z_i M_n^{(i)}> u_n\}&=& o(1), ~~n\to\infty.
\end{eqnarray}
It holds
\begin{eqnarray*}\label{24}
&&P\{M_n^*(z,d) 
\le u_n\}=1-P\{M_n^*(z,d) 
> u_n\}
= 1-P\{z_dM_n^{(d)}> u_n\}\nonumber
\\
&-&\sum_{j=1}^{d-1} P\{z_{j}M_n^{(j)}>u_n, \max_{j+1\le i\le d}(z_{i}M_n^{(i)})\le u_n\}.
\end{eqnarray*}
Note that $\sum_{j=1}^{d-1} P\{z_{j}M_n^{(j)}>u_n, \max_{j+1\le i\le d}(z_{i}M_n^{(i)})\le u_n\}=o(1)$ is equivalent to (\ref{11b}) for $i=d$.
If (\ref{11b}) holds, then we get
\begin{eqnarray*}
&&P\{M_n^*(z,d) 
\le u_n\}= P\{z_{d}M_n^{(d)}\le u_n\}+o(1), ~~n\to\infty.
\end{eqnarray*}
 Then by (\ref{12}) and (\ref{12a}) the expression
\begin{eqnarray*}\label{22ab}
P\{M_n^*(z) \leq u_n\} & = & P\{z_d M_n^{(d)} \leq u_n\}(1+o(1))\nonumber
\\
&=&\exp(-\theta_d (z_d/y)^{k_1})(1+o(1))
\end{eqnarray*}
that is required to obtain the extremal index $\theta_d$ for $Y_n^*(z)$ follows by (\ref{2}) and (\ref{quantile}). 
For some $i\in\{1,...,d\}$ such that (\ref{11b}) holds,  $Y_n^*(z)$ has the extremal index $\theta_i$ for reasons of symmetry. 
The same result 
holds for $Y_n(z)$ due to (\ref{8b}), (\ref{17}) and (\ref{17a}) if the conditions (A2) are additionally fulfilled.
\subsection{Proof of Theorem \ref{T4}} 
\subsubsection{Case (i)}

 \paragraph{Tail index}
 Similarly to the proof of Theorem 4 in Markovich and Rodionov (2020)  we have for all $m\ge 1$
 \begin{eqnarray}\label{47a}&& (z_1/x)^{k_1}\ell_1(x)(1+o(1))\le P\{z_{1}Y_{n,1}>x\}\le P\{Y_{m}^*(z,N_m)>x\} \nonumber
 \\
 &\leq & P\{Y_{m}(z,N_m)>x\}
 \leq
P\{Y_{m}(z,l(x))>x\}+P\{N_m>l(x)\},
\end{eqnarray}
 where $l(x)$ is selected in such a way to neglect the term $P\{N_m>l(x)\}$ taking into account (\ref{15a}). Namely, $l_m$ is replaced by an arbitrary natural number $l=l(x)$ and for all sufficiently large $x$  we set $l(x) = [x^{k_1/\alpha + \delta_1}]$ with $\delta_1>0$.
  \\
  We get
  \begin{eqnarray}\label{40}&& P\{Y_m^*(z,l(x))>x\}=P\{z_{1}Y_{n,1}>x\}+\sum_{i=2}^{l}P\{z_{i}Y_{n,i}>x,\max_{1\le j\le i-1}z_{j}Y_{n,j}\le x\}.\nonumber
  \\
  &&
  \end{eqnarray}
  \paragraph{Condition (A3)} Using notation (\ref{51}) we obtain
  \begin{eqnarray}\label{32a}
  \sum_{i=2}^{l}\xi_i &=&\sum_{s=2}^{\lfloor d_n-1\rfloor}\left(\sum_{i=2}^{s}\xi_i+\sum_{i=s+1}^{l}\xi_i\right)P\{d=s\}
   \end{eqnarray}
   since $d$ and $\{Y_{n,i}\}$ are assumed to be independent and  $1<d<d_n=\min(C,l_n)$ holds. By (\ref{18b})
   we have
   \begin{eqnarray}\label{33}&&
   \sum_{i=s+1}^{l}\xi_i\le \sum_{i=s+1}^{l}P\{z_{i}Y_{n,i}>x\}\le O\left(x^{-k_1(1+\delta)}\right),~x\to\infty,
   \end{eqnarray}
   with $\delta\in(0,k_1\chi_0)$. 
   By (A3) it holds
   \begin{eqnarray}\label{21}
\sum_{i=2}^{s}\xi_i&=&\sum_{i=2}^{s}P\{z_iY_{n,i}>x,\max_{1\le j\le i-1}z_jY_{n,j}\le x\}\nonumber
\\&=&P\{\max_{2\le j\le s}(z_{j}Y_{n,j})>x, z_{1}Y_{n,1}\le x\}=o(P\{Y_{n,1}>x\}) 
\end{eqnarray}
for any $s\in\{2,3,...,\lfloor d_n-1\rfloor\}$ as $x\to\infty$. 
\\
  By (\ref{40})-(\ref{21}) it follows 
   \begin{eqnarray}\label{53a} \!\!\!\!\!\!\!\!\!\!\!\!\!\!\!&&P\{Y^*_{m}(z,l(x))>x\}=P\{z_{1}Y_{n,1}>x\}(1+o(1))=(z_1/x)^{k_1}\ell_1(x)(1+o(1)).
   \end{eqnarray}
  Then by (\ref{47a}) $P\{Y_m^*(z,N_m)>x\}$ has the tail index $k_1$. 
 \paragraph{Condition (A1) (or (A2)) and the independence of $N_m$ and $\{Y_{m,i}\}$}  If, in addition to condition (A1) or (A2),
   $N_m$ and $\{Y_{m,i}\}$ are assumed to be  independent, then $P\{Y_m(z,N_m)>x\}$ and $P\{Y^*_m(z,N_m)>x\}$ have the same tail index $k_1$.
  Really,
  by (\ref{18b}), (\ref{34a}), (\ref{34b}) (or (\ref{13b})) we get for $m\ge 1$
 \begin{eqnarray}\label{32}&& P\{Y_m(z,N_m)>x\}\nonumber
 \\
 &=& \sum_{i=1}^{\infty}P\{Y_m(z,i)>x\}P\{N_m=i\}= 
 x^{-k_1}\Omega(d_m,k_1)(1+o(1)),
 \end{eqnarray}
 where \begin{eqnarray*}&&\Omega(d_m,k_1)=
 \\
 &&\!\!\!\!\!\sum_{i=1}^{\infty}\sum_{s=2}^{\lfloor d_m-1\rfloor}\left(\sum_{j=1}^sz_j^{k_1}\ell_j(x)\Ii\{i>s\}+\sum_{j=1}^iz_j^{k_1}\ell_j(x)\Ii\{i\le s\}\right)P\{d=s\}P\{N_m=i\}.
 \end{eqnarray*}
 The same is valid for $P\{Y^*_m(z,N_m)>x\}$ by the same reasons. 
\paragraph{Extremal index}
Let us find the extremal index of 
$Y_n^*(z,N_n)$ and $Y_n(z,N_n)$. We use formula (31) in Markovich and Rodionov (2020)
\begin{eqnarray}\label{25}P\{z_1M_n^{(1)}>u_n\}&\le & P\{M_{n}^*(z,N_n)>u_n\}\nonumber
\\
&\le & P\{M_{n}^*(z,\tilde{l}_n)>u_n\}+P\{N>\tilde{l}_n\}\nonumber
\\
&\le & P\{M_{n}^*(z,\tilde{l}_n)>u_n\}+ \sum_{i=1}^n P\{N_i>\tilde{l}_n\},
\end{eqnarray}
where  $\tilde{l}_n>l_n$ is selected such that the term $\sum_{i=1}^n P\{N_i>\tilde{l}_n\}\le O(n^{1-\alpha(\chi+\delta)})=o(1)$ is neglected as $n\to\infty$, i.e.
\begin{eqnarray}\label{25a}
\tilde{l}_n=n^{\chi+2\delta}, \qquad \delta=(\chi_0-\chi)/3.
\end{eqnarray}
Similarly to (\ref{40}) 
we get
\begin{eqnarray*}\label{55a}
P\{M_n^*(z,\tilde{l}_n)> u_n\}&=&P\{z_{1}M_n^{(1)}>u_n\}\nonumber
\\
&+&\sum_{i=2}^{\tilde{l}_n}P\{z_{i}M_n^{(i)}>u_n,\max_{1\le j\le i-1}z_{j}M_n^{(j)}\le u_n\}.
\end{eqnarray*}
 Denoting $P\{z_{i}M_n^{(i)}>u_n,\max_{1\le j\le i-1}z_{j}M_n^{(j)}\le u_n\}$ as $\zeta_i$ 
 and splitting the sum $\sum_{i=2}^{\tilde{l}_n}$ as in (\ref{32a}), we obtain
 \begin{eqnarray*}\label{34}&&\sum_{s=2}^{\lfloor d_n-1\rfloor}\sum_{i=2}^{s}\zeta_iP\{d=s\}=o(1)
   \end{eqnarray*}
  by (\ref{11b}) and
 \begin{eqnarray*}&&\sum_{s=2}^{\lfloor d_n-1\rfloor}\sum_{i=s+1}^{\tilde{l}_n}\zeta_iP\{d=s\}=o(1)
\end{eqnarray*}
 by (\ref{12a}) as $n\to\infty$. (\ref{12a}) is valid for $\tilde{l}_n$ since $\chi+2\delta<\chi_0$. Hence, it follows
\begin{eqnarray}\label{55b}&&P\{M_n^*(z,\tilde{l}_n)> u_n\}=P\{z_1M_n^{(1)}>u_n\}(1+o(1)).
\end{eqnarray}
Then by (\ref{25})
\begin{eqnarray*}&&P\{M_{n}^*(z,N_n)>u_n\} 
=P\{z_1M_n^{(1)}>u_n\}(1+o(1))
\end{eqnarray*}
follows and the extremal index of $Y_n^*(z,N_n)$ 
is equal to $\theta_1$ (or for reasons of symmetry to $\theta_i$ for some $i\in\{1,...,d\}$ such that (\ref{11b}) holds).
$Y_n(z,N_n)$ has the same extremal index since $P\{M_n^*(z,N_n)>u_n\}\le P\{M_n(z,N_n)>u_n\}\le  P\{M_{n}(z,\tilde{l}_n)>u_n\}+ \sum_{i=1}^n P\{N_i>\tilde{l}_n\}$ by (\ref{25}) holds, and since, in addition, the condition (A1) or (A2) is  assumed due to (\ref{8b}), (\ref{10}), (\ref{9}), (\ref{17}) and (\ref{17a}).
\subsubsection{Case (ii)}
Let $d=d_0$ a.s. Then the same statements as in Case (i) follow. It is enough to replace $s$ in formulas of Case (i) by $d_0$ with probability one.

\subsection{Proof of Corollary \ref{Cor3}}
\subsubsection{Case $k_1^+\le k_1^-$}
We assume that there is a unique "column" series with a minimum tail index $k_1^+$.
 One can derive 
that
\begin{eqnarray}\label{7a}&& P\{Y_n(z,N_n)>u_n^{+}\}=(z_1^+/y)^{k_1^+}n^{-1}(1+o(1)),\nonumber
\\
&&P\{Y_{n}(z, N_n)>u_n^-\}=(z_1^+/y)^{k_1^+}n^{-k_1^+/k_1^-}(1+o(1))
\end{eqnarray}
as $n\to\infty$. Really, since $Y_n(z,N_n)\le Y_n^+(z^+,N_n^+)$ holds, then 
by (\ref{4ab}), (\ref{quantile}), (\ref{8a}) 
with $d=1$ it holds
\begin{eqnarray}\label{6a}
&&P\{Y_n(z,N_n)>u_n^{+}\}\le P\{Y_n^+(z^+,l_{n^+}) 
>u_n^{+}\}+P\{N_n^+>l_{n^+} 
\}\nonumber
\\
&=&P\{Y_n^+(z^+,l_{n^+}) 
>u_n^{+}\}(1+o(1))=(z_1^+/y)^{k_1^+}n^{-1}(1+o(1)),
\end{eqnarray}
\begin{eqnarray}\label{6aa}
P\{Y_{n}(z, N_n)>u_n^-\}&\le & P\{Y_{n}^+(z^+,l_{n^+}) 
>u_n^-\}(1+o(1))\nonumber
\\
&=&(z_1^+/y)^{k_1^+}n^{-k_1^+/k_1^-}(1+o(1)).
\end{eqnarray}
Let us obtain the lower bounds of $P\{Y_{n}(z, N_n)>u_n^{\pm}\}$.
Since it holds
\begin{eqnarray*}\label{18a}&&P\{-Y_n^-(z^-,N_n^-)>u_n^+\varepsilon\}\nonumber
\\
&=&P\{-Y_n^-(z^-,N_n^-)>u_n^+\varepsilon, N_n^-\le l_{n^-}\} 
+P\{-Y_n^-(z^-,N_n)>u_n^+\varepsilon, N_n^-> l_{n^-}\}\nonumber 
\\
&\le &  
P\{Y_n^-(|z^-|,l_{n^-}) 
>u_n^+\varepsilon\}
+P\{N_n^->l_{n^-}\} 
\end{eqnarray*}
and by (\ref{4ab}), (\ref{quantile}), (\ref{8a}) we get for arbitrary 
$\varepsilon>0$, $d=1$ and $k_1^+\le k_1^-$
\begin{eqnarray}\label{21ab}
&&P\{Y_n(z,N_n)>u_n^+\}=P\{Y_n^+(z^+,N_n^+)+Y_n^-(z^{-},N_n^-)>u_n^+\}\nonumber
\\
&\ge & P\{Y_n^+(z^+,N_n^+)>u_n^+(1+\varepsilon),0<-Y_n^-(z^-,N_n^-)\le u_n^+\varepsilon\}\nonumber
\\
&\ge &P\{Y_n^+(z^+,N_n^+)>u_n^+(1+\varepsilon)\}-P\{-Y_n^-(z^-,N_n^-)>u_n^+\varepsilon\}\nonumber
\\
&\ge & P\{z_1^+Y_{n,1}>u_n^+(1+\varepsilon)\}-P\{Y_n^-(|z^-|, l_{n^-})
>u_n^+\varepsilon\} 
-P\{N_n^->l_{n^-}\} 
\nonumber
\\
&= & \left(\left(\frac{z_1^+}{y(1+\varepsilon)}\right)^{k_1^+}n^{-1}-
\left(\frac{|z_1^-|}{y\varepsilon}\right)^{k_1^-}n^{-k_1^-/k_1^+}\right)(1+o(1))\nonumber
\\
&=& \left(\frac{z_1^+}{y(1+\varepsilon)}\right)^{k_1^+}n^{-1}(1+o(1)),
\end{eqnarray}
\begin{eqnarray}\label{21aa}
P\{Y_n(z,N_n)>u_n^-\}&\ge &
\left(\left(\frac{z_1^+}{y(1+\varepsilon)}\right)^{k_1^+}n^{-k_1^+/k_1^-}-
\left(\frac{|z_1^-|}{y\varepsilon}\right)^{k_1^-}n^{-1}\right)(1+o(1))\nonumber
\\
&=& \left(\frac{z_1^+}{y(1+\varepsilon)}\right)^{k_1^+}n^{-k_1^+/k_1^-}(1+o(1)).
\end{eqnarray}
In the same way, (\ref{7a}) is true for $Y^*_{n}(z, N_n)$.
\paragraph{Tail index}
Let us show  that the 
sequence $Y_{n}(z, N_n)$ in (\ref{3b}) has the tail index $k_1^+$ for $k_1^+\le k_1^-$.
As in the proof of Theorem \ref{T4} let us replace $l_{m^{\pm}}$ 
by an arbitrary natural number $l^{\pm}=l^{\pm}(x)=[x^{k_1^{\pm}/\alpha + \delta_1}]$ with $\delta_1>0$ for all sufficiently large $x$. By Theorem \ref{T2} 
we have
\begin{eqnarray}\label{21ac}
&&P\{Y_m(z,N_m)>x\}\le P\{Y_m^+(z^+,l^+(x))>x\}+P\{N_m^+>l^+(x)\}\nonumber
\\
&=&(z_1^+/x)^{k_1^+}\ell_1(x)(1+o(1))+x^{-k_1^+-\alpha\delta_1}\widetilde{\ell}_m(x^{k_1^+/\alpha + \delta_1})\nonumber
\\
&=&(z_1^+/x)^{k_1^+}\ell_1(x)(1+o(1))
\end{eqnarray}
as $x\to\infty$ for all $m\ge 1$. On the other hand, we get similarly to (\ref{21ab})
\begin{eqnarray}\label{21ad}
&&P\{Y_m(z,N_m)>x\}\nonumber
\\
&\ge & P\{z_1^+Y_{m,1}>x(1+\varepsilon)\}-P\{Y_m^-(|z^-|,l_{m^-}) 
>x\varepsilon\}
-P\{N_m^->l_{m^-}\} 
\nonumber
\\
&=&\left(\left(\frac{z_1^+}{x(1+\varepsilon)}\right)^{k_1^+}-\left(\frac{|z_1^-|}{x\varepsilon}\right)^{k_1^-}\right)\ell_1(x)(1+o(1))
- x^{-k_1^--\alpha\delta_1}\widetilde{\ell}_m(x^{k_1^-/\alpha + \delta_1})\nonumber
\\
&=&(z_1^+/x)^{k_1^+}\ell_1(x)(1+o(1)).
\end{eqnarray}
The same is valid for $Y^*_m(z,N_m)$.
Then the first statement follows.
\paragraph{Extremal index}
By  (\ref{7a}) the extremal index of $Y_{n}(z, N_n)$ and of $Y^*_{n}(z, N_n)$ does not exist for $u_n^-$ due to (\ref{1}).
\\
By (\ref{3b}) and Theorem \ref{T1} it follows that $Y_n^*(z, N_n)$ has the extremal index $\theta^+_1$ corresponding to $k_1^+$ for $u_n^+$. The same holds for $Y_n(z, N_n)$ since by (\ref{17}) and (\ref{25})
 \begin{eqnarray}\label{19}P\{M_n^*(z, N_n)\le u_n\}&=&P\{M_n(z, N_n)\le u_n\}(1+o(1))\end{eqnarray}
follows. Really, similarly to (\ref{17}) we get 
\begin{eqnarray*}0&\le & P\{M_n^*(z, l_n)\le u_n\}-P\{M_n(z, l_n)\le u_n\}\nonumber
\\
&\le &\sum_{j=1}^nP\{Y_j^{*+}(z^+)\le u_n,Y_j^-(z^-)+Y_j^+(z^+)>u_n\}\nonumber
\\
&\le & \sum_{j=1}^nP\{Y_j^{*+}(z^+)\le u_n,Y_j^+(z^+)>u_n\},\label{20}
\end{eqnarray*}
where $Y_n(z)$ and $Y_n^*(z)$ denote $Y_n(z,l_n)$ and $Y_n^*(z,l_n)$.
We obtain
\begin{eqnarray}\label{20a}&&P\{Y_j^{*+}(z^+)\le u_n,Y_j^+(z^+)>u_n\}=P\{Y_j^+(z^+)>u_n\}-P\{Y_j^{*+}(z^+)> u_n\}\nonumber
\\
&\le &P\{Y_n^+(z^+)>u_n\}-P\{z_1^+Y_{n,1}> u_n\}=o(1/n)
\end{eqnarray}
by (\ref{quantile}) and (\ref{8a}) with $d=1$. 
Then  $P\{M_n^*(z, l_n)\le u_n^+\}=P\{M_n(z, l_n)\le u_n^+\}(1+o(1))$ holds and by  (22) in Markovich and Rodionov (2020) $P\{M_n^*(z, l_n)\le u_n^+\}=P\{z_1M_n^{(1)}\le u_n^+\}(1+o(1))$ holds. Then (\ref{19}) follows by (\ref{25}).
\subsubsection{Case $k_1^-\le k_1^+$}
We assume that there is a unique "column" series with the minimum tail index $k_1^-$.
Similarly to (\ref{7a}) 
one can get
\begin{eqnarray}\label{7ab}&& P\{-Y_n(z,N_n)>u_n^{-}\}=(|z_1^-|/y)^{k_1^-}n^{-1}(1+o(1)),\nonumber
\\
&&P\{-Y_{n}(z, N_n)>u_n^+\}=(|z_1^-|/y)^{k_1^-}n^{-k_1^-/k_1^+}(1+o(1)), ~~n\to\infty,
\end{eqnarray}
and by (\ref{21ac}), (\ref{21ad}) for all $m\ge 1$
\begin{eqnarray*}&&P\{-Y_m(z,N_m)>x\}=(|z_1^-|/x)^{k_1^-}\ell_1(x)(1+o(1))
\end{eqnarray*}
as $x\to\infty$. The same is valid for $-Y_{n}^{**}(z, N_n)$.
Then
$-Y_{n}(z, N_n)$ and $-Y_{n}^{**}(z, N_n)$ in (\ref{3b}) have the same tail index $k_1^-$. 
\\
By  (\ref{7ab}) the extremal indices of $-Y_{n}(z, N_n)$ and  $-Y_{n}^{**}(z, N_n)$  do not exist for $u_n^+$ due to (\ref{1}). In the same way as for the case $k_1^+\le k_1^-$ one can derive that  $-Y_{n}(z, N_n)$ and $-Y_{n}^{**}(z, N_n)$ have the same extremal index equal to $\theta^-_1$ for $u_n=u_n^-$.

\subsection{Proof of Corollary \ref{Cor4}}
\subsubsection{Case $k_1^+\le k_1^-$}
By (\ref{3b}) all statements of Theorem \ref{T4} for $Y_{n}^*(z, N_n)$ regarding the tail index and the extremal index  in case of $u_n=u_n^+$ are fulfilled. It remains to show the same for $Y_{n}(z, N_n)$.
\paragraph{Tail index}
One can derive that the tail index of $Y_{n}(z, N_n)$ is given by $k_1^+$ in the same way as in
Theorem \ref{T4}
assuming that $N_m^{\pm}$ and $\{Y^{\pm}_{m,i}\}$ are  independent and (A1) (or (A2)) holds. Really, by (\ref{32}) and similarly to 
(\ref{21ab}) we get
\begin{eqnarray*}
&&P\{Y_m(z,N_m)>x\}\le P\{Y^+_m(z^+,N^+_m)>x\}=x^{-k^+_1}\ell_1(x)\Omega(d^+_m,k^+_1)(1+o(1)),
\end{eqnarray*}
\begin{eqnarray*}
P\{Y_m(z,N_m)>x\}&\ge &P\{Y_m^+(z^+,N_m^+)>x(1+\varepsilon)\}-P\{-Y_m^-(z^-,N_m^-)>x\varepsilon\}
\\
&=&P\{Y_m^+(z^+,N_m^+)>x(1+\varepsilon)\}-P\{Y_m^-(|z^-|,N_m^-)>x\varepsilon\}
\\
&=&(x^{-k^+_1}\ell_1(x)\Omega(d^+_m,k^+_1)-x^{-k^-_1}\ell_1(x)\Omega(d^-_m,k^-_1))(1+o(1))
\\
&=&x^{-k^+_1}\ell_1(x)\Omega(d^+_m,k^+_1)(1+o(1))
\end{eqnarray*}
for any $m\ge 1$.
\\
\paragraph{Extremal index}
Let us find the extremal index of $Y_{n}(z, N_n)$. 
Since $Y_n(z,N_n)\le Y_n^+(z^+,N_n^+)$ holds, then it follows  by (\ref{4ab}),
(\ref{quantile}),
(\ref{8a}) 
and (\ref{6aa}) 
\begin{eqnarray}
&&P\{Y_n(z,N_n)>u_n^{-}\}\le
\left(\frac{z^{*+}}{y(1-\varepsilon)}\right)^{k_1^+}
n^{-k_1^+/k_1^-}(1+o(1)) + o(1/n),\label{26a}
\end{eqnarray}
where $(z^{*+})^{k_1^+}=\sum_{i=1}^{d^+}(z_i^+/\varepsilon_i^*)^{k_1^+}$, $u_n^{\pm}$ are determined as in Corollary \ref{Cor3}.
By 
(\ref{21aa}) we obtain
\begin{eqnarray}
&&P\{Y_n(z,N_n)>u_n^{-}\}\ge \left(\frac{z_1^+}{y}\right)^{k_1^+}n^{-k_1^+/k_1^-}(1+o(1)).\label{26}
\end{eqnarray}
The extremal index of $Y_{n}^*(z, N_n)$ and $Y_{n}(z, N_n)$ does not exist for $u_n^{-}$ due to (\ref{1}), (\ref{26a}) and (\ref{26}).
$Y_{n}(z, N_n)$ has the same extremal index as $Y^*_{n}(z, N_n)$ for $u_n^+$ since (\ref{19}) holds. Really, note that
\begin{eqnarray}&&P\{M_n^*(z, l_n)\le u_n^+\}=P\{M_n(z, l_n)\le u_n^+\}(1+o(1))
\end{eqnarray}
holds
since  (\ref{20a}) is fulfilled by (A1) (see, (\ref{10}), (\ref{9})) or (A2) (see, (\ref{8b})). Then (\ref{19}) follows by (\ref{25}) and (\ref{55b}) in the same way as in Theorem \ref{T4}.
\subsubsection{Case $k_1^-\le k_1^+$}
The proof of the tail index for $-Y_n(z,N_n)$ and $-Y_n^{**}(z,N_n)$ is similar to the previous case. 
The extremal index of $-Y_{n}(z, N_n)$ and  $-Y_{n}^{**}(z, N_n)$  does not exist for $u_n^+$. The extremal index  for  $-Y_{n}(z, N_n)$ and $-Y_{n}^{**}(z, N_n)$ and  $u_n=u_n^-$ is determined symmetrically as for the case $k_1^+\le k_1^-$.
\subsection{Proof of Theorem \ref{T5}}
The proof is similar to that of Theorem \ref{T1}. 
We only indicate the modifications.
\paragraph{Case (i)}
In conditions of Theorem \ref{T1} it holds
\begin{eqnarray}\label{22}P\{Y_{n}(z,l_n)> u_n\} &=& P\{Y_{n}^*(z,l_n)> u_n\}(1+o(1))=P\{z_1Y_{n,1}> u_n\}(1+o(1))\nonumber
\\
&=&(z_1/y)^{k_1} n^{-1}(1+o(1))\end{eqnarray}
as $\ n\to\infty$ (see also (8) and (10) in Markovich and Rodionov (2020)). Hence, (\ref{4b}) implies
\begin{eqnarray*}P\{Y_{n}(z,l_n)> u_n\}&=& o\left(P\{N_n>l_n\}\right).\end{eqnarray*}
Indeed,  $\alpha\chi<1$ follows by (\ref{4b}).
\paragraph{Tail index} Let us find the tail index of $Y^*_{n}(z,N_n)$ and $Y_{n}(z,N_n)$.
The claim
$\alpha\chi_0>1$ is required in Theorem \ref{T1} for the following result derived in formula
(15) in Markovich and Rodionov (2020) 
\begin{eqnarray*}&& P\{Y_{m}(z,l_m)>x\}=P\{Y_{m}^*(z,l_m)>x\}(1+o(1))
\\
&=& P\{Y_{m,1}>x/z_1\}(1+o(1))=(z_1/x)^{k_1}\ell_1(x)(1+o(1))
\end{eqnarray*}
for any fixed $m\ge 1$ and as $x\to\infty$.
To this end, we use 
$l(x)=[x^{k_1/\alpha+\delta_1}]$ with $\delta_1>0$
%
to replace $l_m$ in (16) by Markovich and Rodionov (2020) for sufficiently large $x$ such that $l(x)<x^{k_1\chi_0-\delta_2}$ with arbitrary 
$\delta_2\in(0, k_1\chi_0)$, where $\delta_1+\delta_2=\delta^*$, $\delta^*$ is determined by (\ref{4d}).
Thus, by (\ref{4d}) 
$x^{k_1/\alpha+\delta_1}<x^{k_1\chi_0-\delta_2}$
holds for $x>1$. 
In the same way as in (\ref{13a}) for $d=1$ and by (\ref{11a})
we have
\begin{eqnarray}\label{7c}P\{Y_{m}^*(z,N_m)>x\}&\le &P\{Y_{m}(z,N_m)>x\}\nonumber
\\
&\le & P\{Y_{m}(z,l(x))>x\}+P\{N_m>l(x)\}
\\
&=& (z_1/x)^{k_1}\ell_1(x)(1+o(1))+ x^{-k_1-\alpha \delta_1
}\tilde{\ell}_m(x^{k_1/\alpha +\delta_1})\nonumber 
\\
&=& (z_1/x)^{k_1}\ell_1(x)(1+o(1)),
\nonumber
\end{eqnarray}
\begin{eqnarray}\label{7d}P\{Y_{m}^*(z,N_m)>x\}&\ge & P\{Y_{m,1}>x/z_1\}(1+o(1))\nonumber
\\
&=&(z_1/x)^{k_1}\ell_1(x)(1+o(1))
\end{eqnarray}
as $x\to\infty$. Then the first statement of (i) in Theorem \ref{T5} follows.
\paragraph{Extremal index}
Let us find the extremal index of $Y_{m}(z,N_m)$ and $Y_{m}^*(z,N_m)$. 
Let us use $\widetilde{l}_n
$ as in (\ref{25a}).
By (\ref{27}) 
$\widetilde{l}_n>l_n$ follows and  (\ref{22})  holds for $\widetilde{l}_n$ since $\chi+2\delta<\chi_0$ holds.
By (\ref{15a}), (\ref{27}), (\ref{4b}), (\ref{quantile}), (\ref{22}) 
we have
\begin{eqnarray*}\label{7}
&&(z_1/y)^{k_1} n^{-1}(1+o(1))=P\{z_1Y_{n,1}>u_n\}\le P\{Y^*_n(z,N_n)>u_n\}\nonumber
\\
&\le & P\{Y_n(z,N_n)>u_n\}\le P\{Y_n(z,\widetilde{l}_n)>u_n\}+P\{N_n>\widetilde{l}_n\}\nonumber
\\
&\le &(z_1/y)^{k_1} n^{-1}(1+o(1))+n^{-\alpha(\chi+2\delta)} \tilde{\ell}_n(n^{\chi+2\delta})\nonumber
\\
&=& (z_1/y)^{k_1} n^{-1}(1+o(1)),~~ n\to\infty
\end{eqnarray*}
since
\begin{eqnarray}\label{29}&&\alpha(\chi+2\delta)>1\end{eqnarray} is valid by (\ref{4e}) and $\alpha\chi<1$.
\\
Let us denote $\tilde{l}_n^*=n^{\chi+2.5\delta}$. By (31) in Markovich and Rodionov (2020) we get
\begin{eqnarray}\label{30}P\{z_1M_n^{(1)}>u_n\}&\le & P\{M_{n}^*(z,N_n)>u_n\}
\le   P\{M_{n}^*(z,\tilde{l}_n^*)>u_n\}+P\{N>\tilde{l}_n^*\}\nonumber
\\
 &\le & P\{M_{n}^*(z,\tilde{l}_n^*)>u_n\}+ \sum_{i=1}^n P\{N_i>\tilde{l}_n^*\}.
\end{eqnarray}
We get
\begin{eqnarray}\label{31}
\sum_{i=1}^n P\{N_i>\tilde{l}_n^*\}&=&\sum_{i=1}^nn^{-\alpha(\chi+2\delta)}n^{-0.5\alpha\delta}\tilde{\ell}_i(n^{\chi+2.5\delta})\nonumber
\\
&\le & O(n^{1-\alpha(\chi+2\delta)})=o(1)
\end{eqnarray}
since
 $n^{-0.5\alpha\delta}\tilde{\ell}_i(n^{\chi+2.5\delta})\to 0$ uniformly for all $i\in\{1,...,n\}$ as $n\to\infty$ and (\ref{29}) holds. By  (22) in Markovich and Rodionov (2020), i.e.
 $P\{M_{n}^*(z)\le u_n\}=P\{z_1M_{n}^{(1)}\le u_n\}(1+o(1))$, (\ref{30}) and (\ref{31}) it follows
\begin{eqnarray*}&&\lim_{n\to\infty}P\{M_{n}^*(z,N_n)\le u_n\}=\exp\{-\theta_1(z_1/y)^{k_1}\}
\end{eqnarray*}
and the extremal index of $Y^*_n(z,N_n)$ is equal to $\theta_1$ by (\ref{22}). It remains to show the same for $Y_n(z,N_n)$. To this end, one can derive
\begin{eqnarray*}\label{76a}
\lim_{n\to\infty}P\{M_{n}(z,N_n)\le u_n\}&=&\lim_{n\to\infty}P\{M_{n}^*(z,N_n)\le u_n\}
\end{eqnarray*}
by (\ref{30}), (\ref{31}) and since $P\{M_n^*(z)\le u_n\} = P\{M_n(z)\le u_n\} (1+o(1))$ holds by (23) in Markovich and Rodionov (2020).
\paragraph{Case (ii)} 
The tail index of $P\{Y_{m}(z,N_m)>x\}$ and $P\{Y_{m}^*(z,N_m)>x\}$ is equal to $k_1$ by (\ref{7c}) and (\ref{7d}). Let us find their extremal index. By (\ref{4c}) $\alpha\chi=1$ and then (\ref{31}) hold. Thus, the statements of Item (i) regarding the extremal index are fulfilled.
\subsection{Proof of Theorem \ref{T6}}
\paragraph{Tail index}
If (A3) holds, then $Y_n^*(z,N_n)$ has the tail index $k_1$ by (\ref{53a}) and (\ref{7c}).
As in the proof of Theorem \ref{T4}  the sequences
$Y_n^*(z,N_n)$ and $Y_n(z,N_n)$ have the same tail index $k_1$ if $N_n$ and $\{Y_{n,i}\}$ are  independent and
(A1) (or (A2)) holds by (\ref{32}).
\paragraph{Extremal index}
By (\ref{55b}), (\ref{29}), (\ref{30}) and (\ref{31}) it follows $P\{M^*_{n}(z,N_n)> u_n\}=P\{z_1M^{(1)}_n>u_n\}(1+o(1))$.
 The same is valid for $P\{M_{n}(z,N_n)$ by the same reasoning as in Theorem \ref{T4}. 

\begin{acknowledgements}
 The author was
  supported by the Russian Science Foundation
(grant \mbox{No.\,22-21-00177}). The author would like to thank Igor Rodionov for his useful comments.
\end{acknowledgements}

\end{document}